\documentclass[a4paper,reqno,12pt]{amsart}
\usepackage{amsmath}
\usepackage{amssymb}
\usepackage{amsthm}
\usepackage{mathrsfs}
\usepackage{times}
\allowdisplaybreaks
\allowdisplaybreaks
\theoremstyle{plain}
 \newtheorem{theorem}{Theorem}

 \newtheorem{lemma}[theorem]{Lemma}
\theoremstyle{definition}
 
\theoremstyle{remark}
 
\newcommand{\s}{\smallskip}

\begin{document}
%%%%%%
%%%%%% title, author(s) and etc
%%%%%%
\title[Bilinear Estimates]{ 
Bilinear Estimates\\Associated to the Schr\"odinger Equation\\with a Nonelliptic Principal Part} 
\author {Eiji ONODERA} 
\address{Mathematical Institute, Tohoku University, Sendai 980-8578, Japan}
\email{sa3m09@math.tohoku.ac.jp}
\thanks{The author is supported by 
the JSPS Research Fellowships for Young Scientists 
and the JSPS Grant-in-Aid for Scientific Research 
No.19$\cdot$3304.}
%%%%%%
%%%%%% abstract
%%%%%%
\begin{abstract} 
We discuss bilinear estimates of 
tempered distributions 
in the Fourier restriction spaces 
for the two-dimensional Sch\"odinger equation 
whose principal part is the d'Alembertian.  
We prove that the bilinear estimates hold if and only if 
the tempered distributions are functions.
\end{abstract}
%%%%%%
%%%%%% maketitle
%%%%%%
\maketitle
%%%%%%
%%%%%% Section 1. Introduction
%%%%%%
\section{Introduction}
\label{section:introduction}
This paper is devoted to studying bilinear estimates of 
tempered distributions in the Fourier restriction spaces 
related with the two-dimensional Schr\"odinger equation 
whose principal part is the d'Alembertian. 
The Fourier restriction spaces were originated 
by Bourgain in his celebrated papers 
\cite{BOURGAIN1} and \cite{BOURGAIN2} 
to establish time-local or time-global well-posedness 
of the initial value problem for 
one-dimensional nonlinear Schr\"odinger equations 
and the Korteweg-de Vries equation 
in $L^2(\mathbb{R})$ respectively. 
Generally speaking, 
to solve the initial value problem for 
nonlinear dispersive partial differential equations 
which can be treated by the classical energy method, 
one usually analyzes the interactions of propagation 
of singularities in nonlinearity in detail, 
and applies the regularity properties of free propagators 
to the resolution of singularities. 
It is well-known that propagators of some classes of 
linear dispersive equations with constant coefficients have 
local smoothing effects (see, e.g., \cite{CHIHARA}), 
and dispersion properties 
(see, e.g., \cite{GV} and \cite{STRICHARTZ}). 
Surprisingly, 
the Fourier restriction spaces automatically 
work for both of the analysis of the interactions of 
propagation of singularities 
in the frequency space and 
the application of the regularity properties of 
free propagators. 
For this reason, many applications and refinements of 
the method of the Fourier restriction spaces 
have been investigated in the last decade; 
see, e.g., 
\cite{CDKS}, %\cite{KPV1,KPV2,NTT,TAO1,TAO2}
\cite{KPV1}--\cite{NTT}, \cite{TAO1,TAO2} 
and references therein. 

\s
Here we state the definition of the Fourier restriction spaces. 
The Fourier transform of a function $f(x,t)$ of 
$(x,t)=(x_1,\dotsc,x_n,t)\in\mathbb{R}^{n+1}$ 
is defined by 
$$
\tilde{f}(\xi,\tau)
=
(2\pi)^{-\frac{n+1}{2}}
\iint_{\mathbb{R}^{n+1}}
e^{-it\tau-ix\cdot\xi}
f(x,t)
dxdt, 
$$
where $i=\sqrt{-1}$, 
$(\xi,\tau)=(\xi_1,\dotsc,\xi_n,\tau)\in\mathbb{R}^{n+1}$ 
and 
$x\cdot\xi=x_1\xi_1+\dotsb+x_n\xi_n$. 
Let $a(\xi)$ be a real polynomial of \vspace{-0.05cm}
$\xi=(\xi_1,\dotsc,\xi_n)\in\mathbb{R}^n$. 
Set $\partial_t=\frac{\partial}{\partial{t}}$,
$\partial_j=\frac{\partial}{\partial{x_j}}$, 
$D_t=-i\partial_t$, 
$D_j=-i\partial_j$, 
$D=(D_1,\dotsc,D_n)$, 
$\lvert{\xi}\rvert=\sqrt{\xi\cdot\xi}$, 
$\langle{\tau}\rangle=\sqrt{1+\tau^2}$,  
and 
$\langle\xi\rangle=\sqrt{1+\lvert\xi\rvert^2}$.  
For $s,b\in\mathbb{R}$, the Fourier restriction space 
$X^{s,b}=X^{s,b}(\mathbb{R}^{n+1})$ 
associated to the differential operator 
$D_t-a(D)$ 
is the set of all tempered distributions 
$f$ on $\mathbb{R}^{n+1}$ satisfying 
$$
\lVert{f}\rVert_{s,b}
=
\left(
\iint_{\mathbb{R}^{n+1}}
\big\lvert
\langle{\tau-a(\xi)}\rangle^b
\langle{\xi}\rangle^s
\tilde{f}(\xi,\tau)
\big\rvert^2
\,d\xi
\,d\tau
\right)^{\frac{1}{2}}
<+\infty. 
$$
The free propagator $e^{ita(D)}$ 
of a differential equation $(D_t-a(D))u=0$ 
is defined by 
$$
e^{ita(D)}\phi(x)
=
(2\pi)^{-\frac{n}{2}}
\int_{\mathbb{R}^n}
e^{ix\cdot\xi+ita(\xi)}
\hat{\phi}(\xi)
\,d\xi,
$$
where $\hat{\phi}$ is the Fourier transform of $\phi$ in 
$x\in\mathbb{R}^n$, that is, 
$$
\hat{\phi}(\xi)
=
(2\pi)^{-\frac{n}{2}}
\int_{\mathbb{R}^n}
e^{-ix\cdot\xi}
\phi(x)
\,dx.
$$

\s
In one-dimensional case,  
bilinear estimates in the Fourier restriction spaces 
associated to $D_t-D^2$ and $D_t-D^3$ were completed. 
More precisely, in \cite{KPV1} and \cite{KPV2}, 
Kenig, Ponce and Vega refined the bilinear estimates 
in the Fourier restriction spaces 
with some negative indices $s<0$. 
Nakanishi, Takaoka and Tsutsumi in \cite{NTT} 
constructed sequences of tempered distributions 
breaking the bilinear estimates 
to show the optimality of the indices $s<0$ 
used in \cite{KPV1} and~\cite{KPV2}. 

\s
In \cite{TAO1}
Tao investigated the bilinear estimates associated to 
$a(\xi)=\lvert\xi\rvert^2$ with $n\geqslant2$. 
He dealt with some equivalent estimates 
of the integral of trilinear form, 
and pointed out that the worst singularity occurs 
when an orthogonal relationship 
of three phases in that integral holds. 
Particularly in case $n=2$, 
Colliander, Delort, Kenig and Staffilani 
succeeded in overcoming this difficulty 
by the dyadic decomposition in not only the sizes of phases 
but also the angles among them. 
See \cite{CDKS} for the detail. 
Combining the above results for 
$a(\xi)=\lvert\xi\rvert^2$ with $n=1,2$, 
we have the following. 

\begin{theorem}[\cite{CDKS,KPV1}]
\label{theorem:CDKPSV}
Let $n=1,2$, and let $a(\xi)=\lvert\xi\rvert^2$.
\begin{itemize}
\item[\rm (i)]{
For any $s\in (-\frac{3}{4},0]$, 
there exist $b\in(\frac{1}{2},1)$ and $C>0$ such that 
\begin{align}
  \lVert{uv}\rVert_{s,b-1}
& \leqslant
  C
  \lVert{u}\rVert_{s,b}\,
  \lVert{v}\rVert_{s,b}, 
\label{equation:b1}
\\
  \lVert{\bar{u}\bar{v}}\rVert_{s,b-1}
& \leqslant
  C
  \lVert{u}\rVert_{s,b}\,
  \lVert{v}\rVert_{s,b}. 
\label{equation:b2}
\end{align}
}
\item[\rm (ii)]{
For any $s\in (-\frac{1}{4},0]$, 
there exist $b\in(\frac{1}{2},1)$ and $C>0$ such that 
\begin{equation}
\lVert{\bar{u}v}\rVert_{s,b-1}
\leqslant
C
\lVert{u}\rVert_{s,b}\,
\lVert{v}\rVert_{s,b}. 
\label{equation:b3}
\end{equation}
}
\item[\rm (iii)]{
For any $s<-\frac{3}{4}$ and for any $b\in\mathbb{R}$, 
the estimates 
\eqref{equation:b1} 
and 
\eqref{equation:b2} 
fail to hold, 
and for any $s<-\frac{1}{4}$ and for any $b\in\mathbb{R}$, 
\eqref{equation:b3} 
fails to hold.}
\end{itemize}
\end{theorem}
%
%
%par
Here we mention a few remarks. 
First, the difference between (i) and (ii) 
are basically due to the structure of the products. 
In view of H\"ormander's theorem  
concerned with the microlocal condition 
on the multiplication of distributions 
(see \cite[Theorem~0.4.5]{SOGGE} for instance), 
$u\bar{u}$ needs more smoothness of $u$ 
than $u^2$ and $\bar{u}^2$ 
to make sense. 
Secondly, 
the local smoothing effect and the dispersion property 
of the fundamental solution $e^{it\lvert{D}\rvert^2}$ 
are strongly reflected in these bilinear estimates. 
These are applied to solving the initial value problem 
for some nonlinear Schr\"odinger equations 
in a class of tempered distributions 
which are not necessarily functions. 
Indeed, 
by using the technique developed in \cite{KPV1} 
together with the estimates 
\eqref{equation:b1}, 
\eqref{equation:b2} 
and 
\eqref{equation:b3},  
one can prove time-local well-posedness of 
the initial value problem 
for quadratic nonlinear Schr\"odinger equations of the form 
\begin{alignat}{2}
  D_tu
  -
  \lvert{D}\rvert^2u
& = 
  N_j(u,u)
&
  \quad\text{in}\ 
& \mathbb{R}^n\times\mathbb{R}, 
\label{equation:pde1} 
\\
  u(x,0)
& =
  u_0(x)
&
  \quad\text{in}\ 
& \mathbb{R}^n, 
\label{equation:data1}
\end{alignat}
in Sobolev space $H^s(\mathbb{R}^n)$ 
with $s\in (-\frac{3}{4},0]$ for $j=1,2$ and 
$s\in (-\frac{1}{4},0]$ for $j=3$, respectively. 
Here $n=1,2$, 
$u(x,t)$ is a complex-valued unknown function of $(x,t)$, 
$u_0$ is a given initial data, 
$N_1(u,v)=uv$, 
$N_2(u,v)=\bar{u}\bar{v}$, 
$N_3(u,v)=\bar{u}v$, 
$H^s(\mathbb{R}^n)=\langle{D}\rangle^{-s}L^2(\mathbb{R}^n)$, 
and 
$L^2(\mathbb{R}^n)$ is the set of 
all square-integrable functions on $\mathbb{R}^n$. 

\s
Some two-dimensional nonlinear dispersive equations 
with a nonelliptic principal part 
arise in classical mechanics.  
For example, the Ishimori equation~(\cite{ISHIMORI})
\begin{align*}
D_tu-(D_1^2-D_2^2)u
&=
\frac{-2\bar{u}}{1+\lvert{u}\rvert^2}
\Bigl((D_1u)^2-(D_2u)^2\Bigr)
+
i(D_2\phi{D_1u}+D_1\phi{D_2u}),\\
\phi
&=
-4i\lvert{D}\rvert^{-2}
\left(
\frac{D_1\bar{u}D_2u-D_1uD_2\bar{u}}{1+\lvert{u}\rvert^2}
\right),
\end{align*}
and the hyperbolic--elliptic Davey-Stewartson equation 
(\cite{DS}) \vspace{-0.05cm}
$$
D_tu-(D_1^2-D_2^2)u
=
-\lvert{u}\rvert^2u
-uD_1^2\lvert{D}\rvert^{-2}(\lvert{u}\rvert^2) \vspace{-0.05cm}
$$ 
are well-known two-dimensional nonlinear dispersive equations. 
It is easy to see that 
$e^{it(D_1^2-D_2^2)}$ has 
exactly the same local smoothing and dispersion properties of 
$e^{it(D_1^2+D_2^2)}$ 
since $a(\xi)=\xi_1^2\pm\xi_2^2$ 
are two-dimensional nondegenerate quadratic forms. 
If the gradient $a^\prime(\xi)$ of a quadratic form $a(\xi)$ 
does not vanish for $\xi\ne0$, 
then $e^{ita(D)}$ gains $\frac{1}{2}\,$-spatial differentiation 
globally in time and locally in space. 
If the Hessian $a^{\prime\prime}(\xi)$ of 
an $n$-dimensional quadratic form $a(\xi)$ 
is a nonsingular matrix, 
then the distribution kernel of $e^{ita(D)}$ 
in $\mathbb{R}^n\times\mathbb{R}^n$ is estimated by 
$O(\lvert{t}\rvert^{-\frac{n}{2}})$ for all $t\in\mathbb{R}$
(see, e.q., \cite{KPV0}).  
Then, we expect that 
the bilinear estimates for $a(\xi)=\xi_1^2-\xi_2^2$ 
are the same as those for $a(\xi)=\xi_1^2+\xi_2^2$. 
The purpose of this paper is to examine this expectation. 
However, our answer is negative. 
More precisely, our results are the following.
\begin{theorem}
\label{theorem:main}
Let $n=2$, and let $a(\xi)=\xi_1^2-\xi_2^2$.
\begin{itemize}
\item[\rm (i)]{
For $s \geq 0$, there exists $b\in(\frac{1}{2},1)$ and $C>0$ 
such that the estimates 
\eqref{equation:b1}, 
\eqref{equation:b2} 
and 
\eqref{equation:b3} 
hold.
}
\item[\rm (ii)]{
For any $s<0$ and for any $b\in\mathbb{R}$, 
the estimates 
\eqref{equation:b1}, 
\eqref{equation:b2} 
and 
\eqref{equation:b3} fail to hold.
}
\end{itemize}
\end{theorem}
Note that our results are independent of 
the structure of products. 
In other words, our results depend only on 
the properties of $a(\xi)$, in particular, 
on the noncompactness of the zeros of $a(\xi)$.

\s 
We shall prove Theorem~\ref{theorem:main} in the next section. 
On one hand, 
we directly compute trilinear forms in the phase space 
to show (i) of Theorem~\ref{theorem:main}. 
We see that the Strichartz estimates work for making use of 
the regularity property of the free propagator 
$e^{it(D_1^2-D_2^2)}$ to prove (i). 
\par
On the other hand, 
to prove (ii) of Theorem~\ref{theorem:main},  
we construct two sequences of real-analytic functions 
for which the bilinear estimates break down. 
We observe that one cannot make full use of 
the regularity properties of $e^{it(D_1^2-D_2^2)}$ 
for the negative index $s$. 
More precisely, if $s<0$, then 
these properties cannot work effectively  
near the set of zeros of $a(\xi)$, 
that is, the hyperbola in $\mathbb{R}^2$. 

\s
Finally, we remark that our results seem 
to be strongly related with 
the recent results on bilinear estimates 
of the two dimensional Fourier restriction problems 
by Tao and Vargas in \cite{TV} and \cite{VARGAS}. 
They obtained bilinear estimates of two functions 
restricted on the unit paraboloid in the phase space. 
Their method of proof does not work 
for the restriction on the hyperbolic paraboloid. 
%%%%%%
%%%%%% Section 2. Proof
%%%%%%
\section{Proof of Theorem~\ref{theorem:main}}
\label{section:proof}
Fix $a(\xi)=\xi_1^2-\xi_2^2$. 
Note that $a(\xi)=a(-\xi)$ for any $\xi\in\mathbb{R}^2$. 
First, we prove (i) of Theorem~\ref{theorem:main}. 
Secondly, we prove a lemma needed in the proof of (i). 
Lastly, we conclude this paper by proving 
(ii) of Theorem~\ref{theorem:main}.
%
%\pagebreak
\begin{proof}[Proof of {\rm (i)} of Theorem~{\rm\ref{theorem:main}}]
Let $s \geq 0$ and $\frac{1}{2}<b<1$. 
We employ the idea of trilinear estimates developed in 
\cite{TAO1}.  
In view of the duality argument, 
we have only to show that 
there exists a positive constant $C$ 
depending only on $s$ and $b$ such that 
$$
\lvert{I}\rvert
\leqslant
C
\lVert{f}\rVert_{L^2(\mathbb{R}^3)}\,
\lVert{g}\rVert_{L^2(\mathbb{R}^3)}\,
\lVert{h}\rVert_{L^2(\mathbb{R}^3)},
$$
where 
$$
I
=
\int_{A_1}
\int_{A_2}
\frac{\langle\mu_{0}\rangle^{s}\,
      \langle\mu_{1}\rangle^{-s}\,
      \langle\mu_{2}\rangle^{-s}\,
      f(\mu_0,\tau_0)g(\mu_1,\tau_1)h(\mu_2,\tau_2)}
     {\langle\tau_0+a(\mu_0)\rangle^{1-b}\,
      \langle\tau_1{\pm}a(\mu_1)\rangle^b\,
      \langle\tau_2{\pm}a(\mu_2)\rangle^b}\,
d\tau_0d\tau_1d\tau_2
d\mu_0d\mu_1d\mu_2
$$
and, $A_1$ and $A_2$ are defined by  
\begin{align*}
  A_1
& =
  \{
  (\mu_0,\mu_1,\mu_2)\in\mathbb{R}^{6} 
  \ \vert\  
  \mu_0+\mu_1+\mu_2=0
  \}
\\
  A_2
& =
  \{
  (\tau_0,\tau_1,\tau_2)\in\mathbb{R}^3 
  \ \vert\  
  \tau_0+\tau_1+\tau_2=0 
  \}.
\end{align*}
By using the pairs of signatures 
${\pm}a(\mu_1)$ and ${\pm}a(\mu_2)$ in $I$, 
we can prove 
\eqref{equation:b1}, 
\eqref{equation:b2} 
and 
\eqref{equation:b3} 
together. 
More precisely, 
the pairs $(-,-)$, $(+,+)$ and $(+,-)$ 
correspond to 
\eqref{equation:b1}, 
\eqref{equation:b2} 
and 
\eqref{equation:b3} 
respectively. 
Since 
$\langle\mu_{1}+\mu_{2}\rangle^s
\leqslant
 2^{s}\,
\langle\mu_{1}\rangle^s\,
 \langle\mu_{2}\rangle^s 
$
for $s \geq 0$,
a simple computation gives 
\begin{align*}
  \lvert{I}\rvert
& =
  \biggl\lvert
  \int_{\mathbb{R}^{4}}
  \int_{\mathbb{R}^2}
  \frac{ f(-\mu_1-\mu_2,-\tau_1-\tau_2)
        g(\mu_1,\tau_1)h(\mu_2,\tau_2) }
       { \langle\tau_1{\pm}a(\mu_1)\rangle^b\,
        \langle\tau_2{\pm}a(\mu_2)\rangle^b }
\\ 
&\phantom{=\ } \times
       \frac{ \langle\mu_{1}+\mu_{2}\rangle^{s}\, 
       \langle\mu_{1}\rangle^{-s}\,
       \langle\mu_{2}\rangle^{-s} }
       { \langle-\tau_1-\tau_2+a(\mu_1+\mu_2)\rangle^{1-b} 
}\,
  d\tau_1\,d\tau_2
  \,d\mu_1\,d\mu_2
  \biggr\rvert
\\*[0.2cm]
& \leqslant
  2^s
  \!\int_{\mathbb{R}^{4}}
  \int_{\mathbb{R}^2}
  \frac{\lvert{f(-\mu_1-\mu_2,-\tau_1-\tau_2)}\rvert
        \lvert{g(\mu_1,\tau_1)}\rvert
        \lvert{h(\mu_2,\tau_2)}\rvert}
       {\langle\tau_1{\pm}a(\mu_1)\rangle^b\,
        \langle\tau_2{\pm}a(\mu_2)\rangle^b}\,
  d\tau_1\,d\tau_2
  \,d\mu_1\,d\mu_2
\\*[0.2cm]
& =
  (2\pi)^{-\frac{3}{2}}2^s
  \int_{\mathbb{R}^2}
  \int_{\mathbb{R}}
  \mathscr{F}_{\xi,\tau}^{-1}[\lvert f \rvert](x,t)
  G(x,t)H(x,t)
  \,dt\,dx,
\end{align*}
where
\begin{align*}
  G(x,t)
& =
  \int_{\mathbb{R}^2}
  \int_{\mathbb{R}}
  e^{i(x\cdot\mu+t\tau)}
  \frac{\lvert{g(\mu,\tau)}\rvert}
       {\langle\tau{\pm}a(\mu)\rangle^b}\,
  \,d\tau
  \,d\mu
\\
  H(x,t)
& =
  \int_{\mathbb{R}^2}
  \int_{\mathbb{R}}
  e^{i(x\cdot\mu+t\tau)}
  \frac{\lvert{h(\mu,\tau)}\rvert}
       {\langle\tau{\pm}a(\mu)\rangle^b}\,
  \,d\tau
  \,d\mu,
\end{align*}
and $\mathscr{F}_{\xi,\tau}^{-1}$ 
denotes the inverse Fourier transform on $\xi$ and $\tau$, 
that is, 
$$
\mathscr{F}_{\xi,\tau}^{-1}[\tilde{f}](x,t)
=
(2\pi)^{-\frac{3}{2}}
\iint_{\mathbb{R}^3}
e^{it\tau+ix\cdot\xi}
\tilde{f}(\xi,\tau)
d\xi
d\tau.
$$
The estimates of $G$ and $H$ are the following:
\begin{lemma}
\label{theorem:ste}
For $b>\frac{1}{2}$, 
there exists $C_1=C_1(b)>0$ such that 
for any $g$, $h\in L^2(\mathbb{R}^3)$
$$
\lVert 
G
\rVert_{L^4(\mathbb{R}^3)}
\leqslant
C_1
\lVert 
g
\rVert_{L^2(\mathbb{R}^3)},
\quad
\lVert 
H
\rVert_{L^4(\mathbb{R}^3)}
\leqslant
C_1
\lVert 
h
\rVert_{L^2(\mathbb{R}^3)},
$$
where $L^4(\mathbb{R}^3)$ is the set of 
all Lebesgue measurable functions of 
$(x,t)\in\mathbb{R}^2\times\mathbb{R}$ 
satisfying 
$$
\lVert{F}\rVert_{L^4(\mathbb{R}^3)}
=
\left(
\int_{\mathbb{R}^2}
\int_{\mathbb{R}}
\left\lvert{F(x,t)}\right\rvert^4
dtdx
\right)^{\frac{1}{4}}
<+\infty. 
$$
\end{lemma}
By using Lemma~\ref{theorem:ste}, 
the H\"older inequality and the Plancherel formula,  
we deduce
$$
\lvert{I}\rvert
\leqslant
\lVert{f}\rVert_{L^2(\mathbb{R}^3)}
\lVert{G}\rVert_{L^4(\mathbb{R}^3)}
\lVert{H}\rVert_{L^4(\mathbb{R}^3)}
\leqslant
C
\lVert{f}\rVert_{L^2(\mathbb{R}^3)}
\lVert{g}\rVert_{L^2(\mathbb{R}^3)}
\lVert{h}\rVert_{L^2(\mathbb{R}^3)},
$$
which was to be established.
\end{proof}
\begin{proof}[Proof of Lemma~{\rm\ref{theorem:ste}}]
We show the estimate of $G$. 
Changing a variable by 
$\tau=\lambda{\mp}a(\xi)$, 
we deduce 
\begin{align*}
  G(x,t)
& =
  \int_{\mathbb{R}^2}
  \int_{\mathbb{R}}
  e^{i(x\cdot\xi+t\tau)}
  \frac{\lvert{g(\xi,\tau)}\rvert}
       {\langle\tau{\pm}a(\xi)\rangle^b}
  \,d\tau
  \,d\xi
\\
& =
  \int_{\mathbb{R}}
  e^{it\lambda}
  \langle\lambda\rangle^{-b}
  \left(
  \int_{\mathbb{R}^2}
  e^{ix\cdot\xi}
  e^{{\mp}ita(\xi)}
  \lvert{g(\xi,\lambda{\mp}a(\xi))}\rvert
  \,d\xi
  \right)
  d\lambda
\\
& =
  \int_{\mathbb{R}}
  e^{it\lambda}
  \langle\lambda\rangle^{-b}
  e^{{\mp}ita(D)}
  \psi_{\lambda}(x)
  \,d\lambda,
\end{align*}
where 
$
(\psi_{\lambda})^{\wedge}(\xi)
=
2\pi
\lvert
g(\xi,\lambda \mp a(\xi))
\rvert.
$
Applying the Minkowski inequality, we get
\begin{align}
  \lVert{G}\rVert_{L^4(\mathbb{R}^3)}
& =
  \left(
  \iint_{\mathbb{R}^3}
  \left\lvert
  \int_{\mathbb{R}}
  e^{it\lambda}
  \langle\lambda\rangle^{-b}
  e^{{\mp}ita(D)}
  \psi_{\lambda}(x)
  d\lambda
  \right\rvert^4
  dt\,dx
  \right)^{\!\frac{1}{4}}
\nonumber
\\
& \leqslant
  \int_{\mathbb{R}}
  \left(
  \iint_{\mathbb{R}^3}
  \left\lvert
  e^{it\lambda}
  \langle\lambda\rangle^{-b}
  e^{{\mp}ita(D)}
  \psi_{\lambda}(x)
  \right\rvert^4
  dt\,dx
  \right)^{\!\frac{1}{4}}
  \!d\lambda
\nonumber
\\
& =
  \int_{\mathbb{R}}
  \langle\lambda\rangle^{-b}
  \left(
  \iint_{\mathbb{R}^3}
  \left\lvert
  e^{{\mp}ita(D)}
  \psi_{\lambda}(x)
  \right\rvert^4
  dtdx
  \right)^{\!\frac{1}{4}}
  \!d\lambda.
\label{equation:onnagurui}
\end{align}
\par
Since $a(\xi)$ is a two-dimensional nondegenerate 
quadratic form of $\xi$,
the so-called Strichartz estimate 
$$
\lVert 
e^{\pm ita(D)}u
\rVert_{L^4(\mathbb{R}^3)}
\leqslant
C
\lVert 
u
\rVert_{L^2(\mathbb{R}^2)}
$$
holds (see, e.g., \cite[Appendix]{GS}). 
Using this, 
the Schwarz inequality with $b>\frac{1}{2}$ 
and the Plancherel formula, 
we obtain 
\begin{align*}
  \lVert{G}\rVert_{L^4(\mathbb{R}^3)}
& \leqslant
  C
  \int_{\mathbb{R}}
  \langle\lambda\rangle^{-b}
  \lVert\psi_\lambda\rVert_{L^2(\mathbb{R}^2)}
  \,d\lambda
\\
& \leqslant
  C(b)
  \left(
  \int_{\mathbb{R}}
  \lVert 
  \psi_\lambda
  \rVert_{L^2(\mathbb{R}^2)}^2
  \,d\lambda
  \right)^{\!\frac{1}{2}}
\\
& =
  2\pi
  C(b)
  \left(
  \int_{\mathbb{R}}
  \int_{\mathbb{R}^2}
  \lvert{g(\xi,\lambda{\mp}a(\xi))}\rvert^2
  \,d\xi 
  \,d\lambda
  \right)^{\!\frac{1}{2}}
\\
& =
  2\pi
  C(b)
  \left(
  \int_{\mathbb{R}}
  \int_{\mathbb{R}^2}
  \lvert{g(\xi,\lambda)}\rvert^2
  \,d\xi 
  \,d\lambda
  \right)^{\!\frac{1}{2}}
\\
& =
  2\pi
  C(b)
  \lVert{g}\rVert_{L^2(\mathbb{R}^3)}.
\end{align*}
This completes the proof of Lemma~2.1. 
\end{proof}
\begin{proof}[Proof of {\rm (ii)} of Theorem~{\rm\ref{theorem:main}}] 
Basically we show the optimality in the bilinear estimates 
by constructing suitable Knapp-type counterexamples
as in \cite{KPV1}. 

First, we prove the case $j=1$. 
Fix $s<0$ and $b\in\mathbb{R}$. 
Set $B=\max\{1,\lvert{b}\rvert\}$ for short. 
Suppose that there exists a positive constant $C>0$ 
such that the bilinear estimate \eqref{equation:b1} 
holds for any $u$, $v\in L^2(\mathbb{R}^3)$.  
For $N=1,2,3,\ldots$, set 
$$
\widetilde{u_N}(\xi_1,\xi_2,\tau)
=
\chi_{Q_N}(\xi_1,\xi_2,\tau),
\quad
\widetilde{v_N}(\xi_1,\xi_2,\tau)
=
\chi_{Q_N}(-\xi_1,-\xi_2,-\tau),
$$ 
where $\chi_A$ is the characteristic function of a set $A$, 
and 
$$
Q_N
=
\left\{
(\xi_1,\xi_2,\tau)\in\mathbb{R}^3 
\ \bigg\vert\ 
 N\leqslant \xi_1+\xi_2\leqslant 2N, 
\lvert\xi_1-\xi_2\rvert\leqslant\frac{1}{4N}, 
\lvert\tau\rvert\leqslant\frac{1}{2}
\right\}.
$$
Note that  
\begin{equation}  
Q_N
\subset 
\left\{ 
(\xi_1,\xi_2,\tau)\in \mathbb{R}^3
\ \bigg\vert\ 
\lvert\tau{\pm}a(\xi)\rvert\leqslant1, 
\frac{N}{2}\leqslant\lvert\xi\rvert\leqslant{2N}
\right\},
\label{equation:aya2}
\end{equation}
since 
\begin{gather*}
-\frac{1}{2}
\leqslant
a(\xi)
=(\xi_1+\xi_2)(\xi_1-\xi_2)
\leqslant
\frac{1}{2} \\
\frac{N^2}{2}
\leqslant
\lvert\xi\rvert^2
=
\frac{(\xi_1+\xi_2)^2}{2}
+
\frac{(\xi_1-\xi_2)^2}{2}
\leqslant
2N^2+\frac{1}{32N^2}.
\end{gather*}
By using \eqref{equation:aya2}, 
we deduce 
\begin{align}
  \lVert{u_N}\rVert_{s,b}
& =
  \left(
  \iint_{Q_N}
  \langle\tau-a(\xi)\rangle^{2b}\,
  \langle\xi\rangle^{2s}
  d\tau
  d\xi
  \right)^{\frac{1}{2}}
\nonumber
\\
& \leqslant
  2^{B-s}N^s
  \left(
  \iint_{Q_{N}}d\tau d\xi
  \right)^{\frac{1}{2}}
  \nonumber
\\
& =
  2^{B-s-1}N^s,
\label{equation:ine1}
%\\*[0.3cm]
\intertext{and}
  \lVert{v_N}\rVert_{s,b}
&  =
  \left(
  \iint_{Q_N}
  \langle\tau+a(\xi)\rangle^{2b}\,
  \langle\xi\rangle^{2s}
  d\tau
  d\xi
  \right)^{\frac{1}{2}}
\nonumber
\\
& \leqslant
  2^{B-s}N^s
  \left(
  \iint_{Q_N}d\tau d\xi
  \right)^{\frac{1}{2}}
\nonumber
\\
& =
  2^{B-s-1}N^s.
\label{equation:ine2}
\end{align}
A simple computation shows that 
for large $N\in\mathbb{N}$, 
\begin{align*}
  \widetilde{u_Nv_N}(\xi,\tau)
& =
  (2\pi)^{-\frac{3}{2}}
  \iint_{\mathbb{R}^3}
  \chi_{Q_N}(\xi-\eta,\tau-\lambda)
  \chi_{Q_N}(-\eta,-\lambda)
  d\eta
  d\lambda
\\
& =
  (2\pi)^{-\frac{3}{2}}
  \iint_{Q_N}
  \chi_{Q_N}(\xi+\eta,\tau+\lambda)
  d\eta
  d\lambda
\\
& \geqslant
  \frac{1}{2^{6+\frac{1}{2}}\,\pi^{\frac{3}{2}}}
  \chi_{R_N}(\xi,\tau),
\end{align*}
where
$$
R_N
=
\left\{
(\xi_1,\xi_2,\tau)\in \mathbb{R}^3 
\ \bigg\vert\ 
\lvert\xi_1+\xi_2\rvert\leqslant\frac{N}{2}, 
\lvert\xi_1-\xi_2\rvert\leqslant\frac{1}{8N}, 
\lvert\tau\rvert\leqslant\frac{1}{4}
\right\}.
$$
Since 
$
R_N
\subset
\big\{
(\xi_1,\xi_2,\tau)\in\mathbb{R}^3
\ \big\vert\ 
\lvert\tau{\pm}a(\xi)\rvert\leqslant1, 
\lvert\xi\rvert\leqslant{\frac{N}{2}}
\big\},
$
we get
\begin{align}
  \lVert{u_Nv_N}\rVert_{s,b-1}
& =
  \left(
  \iint_{\mathbb{R}^3}
  \lvert\widetilde{u_{N}v_{N}}(\xi,\tau)\rvert^2
  \langle\tau-a(\xi)\rangle^{2(b-1)}\,
  \langle\xi\rangle^{2s}
  d\tau
  d\xi
  \right)^{\frac{1}{2}}
\nonumber
\\
& \geqslant
  \frac{1}{2^{6+\frac{1}{2}}\pi^{\frac{3}{2}}}
  \left(
  \iint_{R_N}
  \langle\tau-a(\xi)\rangle^{2(b-1)}\,
  \langle\xi\rangle^{2s}
  d\tau
  d\xi
  \right)^{\frac{1}{2}}
\nonumber
\\
& \geqslant
  2^{B-6-\frac{1}{2}}\pi^{-\frac{3}{2}}
  N^s
  \left(
  \iint_{R_N}
  d\tau
  d\xi
  \right)^{\frac{1}{2}}
\nonumber
\\
& =
  2^{B-8-\frac{1}{2}}\pi^{-\frac{3}{2}}
  N^s.
\label{equation:ine3}
\end{align}
Substitute 
\eqref{equation:ine1},
\!\eqref{equation:ine2} 
and 
\eqref{equation:ine3} 
into 
\eqref{equation:b1}. 
\!Then we have 
$
2^{B-8-\frac{1}{2}}\pi^{-\frac{3}{2}}
N^s
\!\leqslant
2^{2B-2s-2}N^{2s}\!, 
$
which becomes 
$
2^{-B+2s-6-\frac{1}{2}}\pi^{-\frac{3}{2}}
\leqslant
N^s.
$
Since $s<0$, 
the right hand side of the above goes to zero 
as $N\rightarrow\infty$ 
while the left hand side is 
a strictly positive constant 
depending only on $s<0$ and $b\in\mathbb{R}$. 
This is contradiction. 
Then, this completes the proof of the case $j=1$.  

\s
The cases $j=2,3$ are proved in the same way.  
Let $Q_N$ be the same as above. 
For $j=2$, set 
$$
\widetilde{u_N}(\xi,\tau)=\chi_{Q_N}(-\xi,-\tau),
\quad
\widetilde{v_N}(\xi,\tau)=\chi_{Q_N}(\xi,\tau), 
$$
and for $j=3$, set 
$$
\widetilde{u_N}(\xi,\tau)=\chi_{Q_N}(-\xi,-\tau),
\quad
\widetilde{v_N}(\xi,\tau)=\chi_{Q_N}(-\xi,-\tau). 
$$
We omit the detail about the cases $j=2,3$. 
\end{proof}
\flushleft{\bf Acknowledgment.} 
The author would like to thank   
Hiroyuki Chihara 
for a number of his valuable suggestions and 
encouragements.

%%%%%%
%%%%%% references
%%%%%%
%\bibliographystyle{amsplain}

\end{document}